\documentclass[12pt]{article}

\usepackage[english]{babel}
\usepackage{amssymb,amsmath}
\usepackage{hyperref}

\newcommand{\footremember}[2]{%
    \footnote{#2}
    \newcounter{#1}
    \setcounter{#1}{\value{footnote}}%
}

\newtheorem{theorem}{Theorem }[section]
\newtheorem{lemma}[theorem]{Lemma}

\newtheorem{proposition}[theorem]{Proposition}

\newcommand{\qed}{\hspace*{\fill}$\Box$}

\title{New necessary conditions for Paley type partial difference sets in Abelian groups}

\author{Zeying Wang \footremember{MTU}{Michigan Technological University}  \footnote{zeying@mtu.edu}}

\date{}

\begin{document}

\maketitle 

\abstract{In this paper we prove that  if there is a regular Paley type partial difference set in an Abelian group $G$ of order $v$, where $v=p_1^{2k_1}p_2^{2k_2}\cdots p_n^{2k_n}$, $n\ge 2$,  $p_1$, $p_2$, $\cdots$, $p_n$ are distinct odd prime numbers, then for any $1 \le i \le n$, $p_i$ is congruent to 3 modulo 4 whenever $k_i$ is odd. These new necessary conditions further limit the specific order of an Abelian group $G$ in which there can exist a Paley type partial difference set. Our result is similar to a result on Abelian Hadamard (Menon) difference sets proved by Ray-Chaudhuri and Xiang in 1997.
 }

\section{Introduction and the main result}

Let $G$ be a finite Abelian group of order $v$, and let $D\subseteq G$ be a subset of size $k$. We say $D$ is a $(v,k,\lambda,\mu)$-{\it partial difference set} (PDS) in $G$ if the expressions $gh^{-1}$, $g$, $h\in D$, $g\neq h$, represent each non-identity element in $D$ exactly $\lambda$ times, and each non-identity element of $G$ not in $D$ exactly $\mu$ times. If we further assume that $D^{(-1)}=D$  (where $D^{(s)}=\{g^s:g\in D\}$) and  $e \notin D$ (where $e$ is the identity element of $G$), then $D$ is called a {\it regular} partial difference set.  A regular PDS is called {\it trivial} if $D\cup\{e\}$ or $G\setminus D$ is a subgroup of $G$. The condition that $D$ be regular is not a very restrictive one, as  $D^{(-1)}=D$ is automatically fulfilled whenever $\lambda\neq\mu$, and $D$ is a PDS if and only if $D\cup\{e\}$ is a PDS. The {\it Cayley graph over $G$ with connection set $D$}, denoted by Cay($G$, $D$), is the graph with the elements of $G$ as vertices, and in which two vertices $g$ and $h$ are adjacent if and only if $gh^{-1}$ belongs to $D$.  When the connection set $D$ is a regular partial difference set, Cay($G$, $D$) is a strongly regular graph.  The importance of regular PDSs lies in the fact that they are equivalent to strongly regular Cayley graphs. For more information on partial difference sets, we refer the reader to a survey of Ma \cite{MA94b}. \\

\medskip

Throughout this paper, we will use the following standard notations: $\beta=\lambda-\mu$ and  $\Delta=\beta^2+4(k-\mu)$.

\medskip

Partial difference sets with parameters $(v, (v-1)/2, (v-5)/4, (v-1)/4)$  are called {\it Paley type partial difference sets}.  Over the last three decades this subject has seen active research, see for example \cite{Arasu}, \cite{Davis}, \cite{Xiang12} (Theorem 3.2),  \cite{MA95}, \cite{Polhill}, \cite{Xiang96}. There are two key problems on Paley type PDSs in Abelian groups:
\begin{itemize}
\item[1.] For what order of the group, can we find Paley-type PDSs?

\item[2.] In which type of groups of given order, can we find Paley-type PDSs?

\end{itemize}
The second question seems largely out of reach , although there is some literature on this topic, see for example \cite{Davis}, \cite{MA95}. In this paper, we will focus on question 1.  It is well-known that when $q \equiv 1 \pmod 4$ and $q$ is a prime power,  the non-zero squares of a finite field ${\mathbb F}_q$ form a  Paley type PDS in the additive group of $\mathbb{F}_q$.  A further important result  was proved by S.L. Ma in 1984:

\begin{theorem}{\rm \cite{MA84}} \label{Ma84}
 Let $D$ be an Abelian regular $(v, k, \lambda, \mu)$-PDS, and assume that $\Delta$ is not a perfect square. Then $D$ is of Paley type; more precisely, $D$ has parameters
$$\left(p^{2s+1}, \frac{p^{2s+1}-1}{2}, \frac{p^{2s+1}-5}{4}, \frac{p^{2s+1}-1}{4}\right),$$
where $p$ is a prime congrent to 1 modulo 4.
\end{theorem}

Let  $D$ be a regular Paley type PDS in an Abelian group $G$, where $|G|=v$.  Then $\Delta=(-1)^2+4(\frac{v-1}{4})=v$.  If $v$ is not a square,  by Theorem \ref{Ma84}, $|G|=v=p^{2s+1}$ for some prime $p\equiv 1 \pmod 4$. For a  prime power $q\equiv 1 \pmod 4$, we can always construct a Paley type PDS in $(\mathbb{F}_q, +)$ using the non-zero squares of the finite field $\mathbb{F}_q$. Thus to answer question 1, we only need to focus on the existence of Paley type PDSs   when $v$ is a perfect square and not a prime power, that is, when $v=p_1^{2k_1}p_2^{2k_2}\cdots p_n^{2k_n}$, $n \ge 2$.

When $|G|=p_1^{2k_1}p_2^{2k_2}\cdots p_n^{2k_n}$, $p_1$, $p_2$, $\cdots$, $p_n$ are distinct odd prime numbers, and  all $k_i$s are even,    Polhill  (\cite{Polhill}) constructed Paley type PDSs in $G=\mathbb{Z}_{p_1}^{2k_1}\times \mathbb{Z}_{p_2}^{2k_2}\times\cdots \times \mathbb{Z}_{p_n}^{2k_n}$.  In this paper we focus on which of the $k_i$s would possibly be odd.

\medskip
We now state our main theorem.

\begin{theorem}\label{Main}
Let $G$ be an Abelian group of order $v$, where $v=p_1^{2k_1}p_2^{2k_2}\cdots p_n^{2k_n}$, $n \ge 2$,  $p_1$, $p_2$, $\cdots$, $p_n$ are distinct odd prime numbers.  If there is a  regular Paley type PDS in $G$, then for any  $1 \le i \le n$, if $k_i$ is odd, we have $p_i \equiv 3 \pmod 4$.
\end{theorem}

At this point, we want to point out a strikingly similar result on Abelian Hadamard (Menon) difference sets proved by  Ray-Chaudhuri and Xiang in 1997.

When $\lambda=\mu$, then a $(v,k,\lambda,\mu)$ partial difference set is called a $(v,k,\lambda)$ difference set (DS). Hadamard (Menon) difference sets, having parameters $(4m^2, 2m^2-m, m^2-m)$, are of particular interest due  to their connections with Hadamard matrices. 

\begin{theorem}{\rm \cite{Xiang97}}\label{Xiang}
 If there is a Hadamard difference set in an Abelian group $G=\mathbb{Z}_2 \times \mathbb{Z}_2 \times P$, where $|P|=p^{2\alpha}$, $\alpha$ is odd, $p$ is an odd prime number, then $p$ is a prime congruent to 3 modulo 4.
\end{theorem}

Although we expect there should be a deeper reason explaining this similarity, we are currently not aware of any general argument.
\section{Proof of the Main Result}

Below we cite three results on Abelian regular partial difference sets. The first of these was proved by S.L. Ma in \cite{MA94a}, the second one was proved by K.T. Arasu, D. Jungnickel, S.L. Ma and A. Pott in \cite{Arasu}, and the last one, the local multiplier theorem, was proved by S. De Winter, E. Kamischke, and Z. Wang in \cite{SDWEKZW}.

\begin{proposition}{\rm \cite{MA94a}}\label{Ma1}
Let D be a nontrivial regular $(v,\,k,\,\lambda,\,\mu)$-PDS in an Abelian group G. Suppose $\Delta$ is a  perfect square.  If $N$ is a subgroup of $G$ such that $\gcd(\left|N\right|, \left|G\right|/\left|N\right|)=1$ and $\left|G\right|/\left|N\right|$ is odd, then $D_1=D\cap N$ is a (not necessarily non-trivial)  regular $(v_1,\,k_1,\,\lambda_1,\,\mu_1)$-PDS with
$$v_1=|N|, \; \beta_1=\lambda_1-\mu_1=\beta-2\theta \pi, \; \Delta_1=\beta_1^2+4(k_1-\mu_1)=\pi^2$$
and
 $$k_1=\frac{1}{2}\left[ |N| +\beta_1\pm \sqrt{(|N|+\beta_1)^2-(\Delta_1-\beta_1^2)(|N|-1)}  \right].$$
where $\pi=\gcd(|N|, \sqrt{\Delta})$ and $\theta$ is the integer satisfying $(2\theta-1)\pi\leq\beta<(2\theta+1)\pi$.
\end{proposition}

\begin{theorem}{\rm \cite{Arasu}}\label{Arasu}
Let $\Gamma$ be a strongly regular Cayley graph based on an Abelian group $G$, with parameters $v$, $k$, $\lambda$, and $\mu$ satisfying $\beta=\lambda-\mu=-1$. Then, up to complementation, $\Gamma$ is either of Paley type or it has parameters $(243, 22, 1, 2)$.
\end{theorem}

\begin{theorem}{\rm \cite{SDWEKZW}}\label{lmt}
Let D be a regular $(v,k,\lambda,\mu)$-PDS in an Abelian group $G$. Furthermore assume $\Delta$ is a perfect square.  Then $g\in G$ belongs to $D$ if and only if $g^s$ belongs to $D$ for all $s$ coprime with $o(g)$, the order of $g$.
\end{theorem}

Using Theorem \ref{Arasu}, we can now prove the following:

\begin{lemma} \label{lemma_inclusion}
 Let $D$ be a  regular Paley type PDS in an Abelian group $G$, where $|G|=v$  is a perfect square.  If $N$ is a
non-trivial subgroup of $G$ such that $\gcd(\left|N\right|, \left|G\right|/\left|N\right|)=1$,  $\left|G\right|/\left|N\right|$ is odd, and $|N|\neq 243$, then $D_1=D\cap N$ is a regular Paley type PDS in $N$.
\end{lemma}

{\bf Proof:} Clearly $\Delta=(\lambda-\mu)^2+4(k-\mu)=1+(v-1)=v$. Applying Proposition \ref{Ma1}, we have  $(2\theta-1)\pi\leq\beta=-1<(2\theta+1)\pi$. Since $\pi=\gcd(|N|, \sqrt{\Delta}) \ge 1$, it follows that $\theta=0$.  Hence $\beta_1=\beta-2\theta \pi=-1$. By Theorem \ref{Arasu},  $D_1=D\cap N$ is a regular Paley type PDS in $N$.\qed\

\begin{lemma}\label{lemma_imp}
Let $G=N\times H$ be an Abelian group, where $N$ and $H$ are subgroups of order $q$ and $p^k$ respectively. Also, assume that $\gcd(p,q)=1$, $p$ is an odd prime number, and $q$ is a positive odd integer.  Let $g=nh\in D$ with $n \in N$ and $h\in H\setminus\{1_H\}$, $o(h)=p^r$.  Then $nh\in D$ if and only if $nh^x \in D$ for all $x$ satisfying $1 \le x \le p^r-1$, and $\gcd(x, p)=1$.
\end{lemma}

{\bf Proof:} Let $x$ be a positive integer satisfying $1 \le x \le p^r-1$ and $\gcd(x,p)=1$. Since $\gcd(p,q)=1$, it follows that $x$, $x+p^r$, $\cdots$, $x+(q-1)p^r$ are in different residue classes modulo  $q$. Thus there exists an integer $t$, $0 \le t \le q-1$, such that $x+tp^r \equiv 1 \pmod q$. Clearly, $\gcd(x+tp^r, q)=1$ and $\gcd(x+tp^r,p^r)=\gcd(x,p^r)=1$. Since $\gcd(p,q)=1$, it follows that  $\gcd(x+tp^r,p^r q)=1$.  As $g=nh \in D$ and ${\rm o}(g)\,|\,p^r q$,  by Theorem \ref{lmt} (the Local Multiplier Theorem), we have 
$$(nh)^{x+tp^r}=nh^x \in D.$$
\qed

\medskip

Now we are ready to prove our main theorem.

\medskip

{\bf Proof:} Here $v=p_1^{2k_1}p_2^{2k_2}\cdots p_n^{2k_n}$, and it is easy to check that $\Delta=v$.

For any $i$, $1 \le i \le n$, we let $G=N \times H$, where $N$ and $H$ are subgroups of $G$ with orders $v/p_i^{2k_i}$ and $p_i^{2k_i}$ repectively. Next we assume that $D$ is a regular Paley type PDS in $G$. As $|N|\neq 3^5=243$, by Lemma \ref{lemma_inclusion}, $D_1=D\cap N$ is a Paley type PDS in $N$ with $\mu_1=\frac{|N|-1}{4}$. Let $n$ be any non-identity element of $N$ and $n \notin D_1$. In total, there are exactly two mutually exclusive types of representations of $n$ as differences from $D$, and these are as follows:
\begin{itemize}
\item[(i)] there are $\mu_1$ representations of $n$ of the form $n_1 n_2^{-1}$ with $n_1, \,n_2\in D_1$;

\item[(ii)]  $n$ can also be written as $(n_1h)(n_2h)^{-1}$ if $n=n_1 n_2^{-1}$, $n_1, \, n_2 \in N$,  and $n_1h \in D$, $n_2h \in D$, $h \in H\setminus\{1_H\}$.

\end{itemize}
By Lemma \ref{lemma_imp}, if $n_1h , \, n_2 h\in D$, so are  $n_1 h^x, \, n_2 h^x$ for any $x$ with $1 \le x \le {\rm o}(h)-1$ and $\gcd(x, {\rm o}(h))=1$. Thus the second type of representations always appear in a set of size $\phi(o(h))$, where $o(h)=p_i^r$ for some $r$ with $1 \le r \le 2k_i$.  As $\phi(o(h))=p_i^r-p_i^{r-1}$ is always divisible by $p_i-1$, we have 
$$\mu=\mu_1+s(p_i-1)=\frac{|N|-1}{4}+s(p_i-1) \;\;\mbox{ for some integer $s$}. $$
On the other hand, $\mu=\frac{|G|-1}{4}$. It follows that
$$\mu-\mu_1=s(p_i-1)=\frac{|G|-|N|}{4}=p_1^{2k_1}p_2^{2k_2}\cdots p_{i-1}^{2k_{i-1}}p_{i+1}^{2k_{i+1}}\cdots p_n^{2k_n}\frac{p_i^{2k_i}-1}{4}.$$
It follows that $p_i-1\,|\,p_1^{2k_1}p_2^{2k_2}\cdots p_{i-1}^{2k_{i-1}}p_{i+1}^{2k_{i+1}}\cdots p_n^{2k_n}\frac{p_i^{2k_i}-1}{4}$.  

Since 
$$p_i^{2k_i}-1=(p_i-1)(p_i^{2k_i-1}+p_i^{2k_i-2}+\cdots+p_i+1),$$
and $p_1$, $p_2$, $\cdots$, $p_n$ are odd prime numbers, 
$$p_i-1\,|\,p_1^{2k_1}p_2^{2k_2}\cdots p_{i-1}^{2k_{i-1}}p_{i+1}^{2k_{i+1}}\cdots p_n^{2k_n}\frac{p_i^{2k_i}-1}{4}$$ holds only when
 \begin{equation}\label{congru_1}
p_i^{2k_i-1}+p_i^{2k_i-2}+\cdots+p_i+1 \equiv 0 \pmod 4.
\end{equation}

From Congruence (\ref{congru_1}), it easily follows that if $k_i$ is odd, we have $p_i \equiv 3 \pmod 4$. This proves the theorem. \qed

\bigskip

\medskip

{\bf Sample Application: } By Theorem \ref{Main}, there does not exist a regular (225, 112, 55, 56)-PDS in Abelian groups since $225=3^2 \times 5^2$ and $5 \equiv 1 \pmod 4$. But there are strongly regular graphs with parameters (225, 112, 55, 56). More generally, there does not exist regular Paley type PDSs of order $25q^2$ in Abelian groups, where $\gcd(5, q)=1$.\\

Given our main result a first natural question is:  Does there exist a regular Paley type PDS  in an Abelian group of order $v=p^{2k}q^2$, where $\gcd(p, q)=1$,  $p$ is a prime number  congruent to 3 modulo 4, and $k$ is odd?    For example, does there exist a regular Paley type PDS of order $3^2 7^2$, or a regular Paley type PDS of order $3^4 7^2$?

\end{document}